\documentclass[a4paper,12pt]{amsart}
\usepackage{amsmath,latexsym,amscd}

\title{Minimal threefolds of small slope and the Noether inequality for
canonically polarized threefolds}
\author{Meng Chen}
\address{Institute of Mathematics, Fudan University,
Shanghai, 200433, PR China} \email{mchen@fudan.edu.cn}
\thanks{This paper is supported by The Institute of Mathematical
Sciences, The Chinese University of Hong Kong. The project is also
supported by the National Natural Science Foundation of China
(No.10131010), Shanghai Scientific $\&$ Technical Commission
(Grant 01QA14042)}
\newcommand{\bQ}{{\mathbb Q}}
\newcommand{\bP}{{\mathbb P}}
\newcommand{\roundup}[1]{\ulcorner{#1}\urcorner}
\newcommand{\rounddown}[1]{\llcorner{#1}\lrcorner}

\newtheorem{thm}{Theorem}[section]
\newtheorem{lem}[thm]{Lemma}
\newtheorem{cor}[thm]{Corollary}
\newtheorem{prop}[thm]{Proposition}
\newtheorem{claim}[thm]{Claim}
\theoremstyle{definition}

\newtheorem{setup}[thm]{}

\theoremstyle{remark}

\begin{document}
\begin{abstract}
Assume that $X$ is a smooth projective 3-fold with ample $K_X$. We
study a problem of Miles Reid to prove the inequality
$$K_X^3\ge \frac{2}{3}(2p_g(X)-5),$$
where $p_g(X)$ is the geometric genus. This inequality is sharp
according to known examples of M. Kobayashi. We also birationally
classify arbitrary minimal 3-folds of general type with small
slope.
\end{abstract}
%%%%%%%%%%%%%%%%%%%%
\maketitle
%%%%%%%%%
\pagestyle{myheadings} \markboth{\hfill M. Chen\hfill}{\hfill
Minimal 3-folds of small slope and the Noether inequality \hfill}
%%%%%%%%%%%%%%%%%%%%%%%%%%%%%%%%%%
\section{\bf Introduction}

We work over an algebraically closed field $k$ of characteristic
0.

On an irreducible complete curve $C$ (may be singular), one has
$\deg(K_C)\ge 2g(C)-2$ where $g(C)$ is the geometric genus of $C$.

This inequality has a 2-dimensional analogue which is the famous
"Noether's inequality" (see \cite{N}). Explicitly, on a minimal
surface $S$ (with RDP singularities) of general type, one has
$$K_S^2\ge 2p_g(S)-4$$
where $p_g(S):=h^0(S, K_S)$ is the geometric genus of $S$.
Together with the Bogomolov-Miyaoka-Yau inequality (cf. \cite{M},
\cite{Y}): $K_S^2\leq 9\chi(O_S)$, they have ever played very
important roles in surface theory (for instance the surface
geography: see \cite{Cat}, \cite{C1}, \cite{C2}, \cite{Per} and
\cite{X} etc.).

The importance of the Noether inequality in mind, Miles Reid first
asked the question seeking for a 3-dimensional analogue in early
1980's. Since then, there have been many papers which give
effective Noether type of inequalities either in the form $K^3\ge
a\chi+b$ (see \cite{O} and  \cite{Ba2} etc.) or for restricted
objects (see \cite{Kob} and \cite{CM} etc.). We mention here that
an effective linear inequality in terms of $\chi$ seems to be
impossible because $\chi$ could be both positive and negative for
a general 3-fold of general type. It is Kobayashi's interesting
examples (\cite{Kob}) that shows the naive inequality $K^3\ge
2p_g-6$ (in 3-dimensional case) is not correct in general. Thus it
becomes more interesting  what the 3-dimensional Noether
inequality is even under very restricted conditions. Such a
"Noether inequality" is by all means important to the 3-fold
geography (see \cite{H}).

The aim of this paper is to present a Noether inequality for
canonically polarized threefolds and to give a classification to
arbitrary minimal 3-folds of general type with small slope. Our
main results are as follows.

\begin{thm}\label{T1} Let $X$ be a smooth projective threefold
with ample canonical divisor $K_X$. Then
$$K_X^3\ge \frac{2}{3}(2p_g(X)-5).$$
\end{thm}

\begin{thm}\label{T1.5} Let $X$ be a minimal projective Gorenstein
3-fold of general type with canonical singularities. Assume
$$K_X^3<2p_g(X)-6.$$
Then $X$ is fibred by curves of genus $2$.
\end{thm}

\begin{thm}\label{T2} Let $X$ be a projective minimal 3-fold with
only canonical singularities. Assume
$$K_X^3<\frac{1}{2}(3p_g(X)-5)$$
(and $p_g(X)\not\in [2,11]$). Then $X$ must be fibred by curves of
genus 2.
\end{thm}

Further classifications to restricted minimal 3-folds are
presented as an interesting application.

\begin{cor}\label{cor} Let $X$ be a smooth projective 3-fold with
ample $K_X$. Assume
$$K_X^3<\frac{3}{2}p_g(X)-\frac{9}{2}.$$
Then $X$ must be canonically fibred by curves of genus 2.
\end{cor}

Theorem \ref{T1} is sharp according to M. Kobayashi's interesting
examples (\cite{Kob}) which say that there are canonically
polarized 3-folds with infinite number of configurations of
invariants $(K^3, p_g)$ satisfying the equality:
$K^3=\frac{2}{3}(2p_g-5)$. Again due to Kobayashi's examples,
Theorem \ref{T1.5} and Theorem \ref{T2} are not empty and  they
are parallel to surface case (see \cite{Reid} and \cite{Xiao}). We
do not know, however, whether both are optimal.

Based on our previous paper (\cite{CM}), and in order to prove
Theorem \ref{T1}, we need to treat the most difficult case, i.e.
when $X$ is canonically fibred by surfaces of general type with
$(c_1^2, p_g)=(1,2)$, through new methods. The new observation of
this paper is that we may choose a special embedded resolution to
the given polarized 3-fold and then successfully apply the
Kawamata-Viehweg vanishing theorem to $\bQ$-divisors on both
3-folds and surfaces to estimate the dimension of those
cohomological groups in question. We bound the $K_X^3$ from below
by studying the bicanonical system rather than in the traditional
way.
\bigskip

This note was written while I was visiting the Institute of
Mathematical Sciences, the Chinese University of Hong Kong. I
would like to thank Eckart Viehweg for his patient explaining my
frequent email queries. Thanks are also due to Kang Zuo for
effective discussions and for his hospitality. Finally I
appreciate many helps from both Keiji Oguiso and Seunghun Lee.

\section{\bf Proof of Theorem \ref{T1} and Theorem \ref{T1.5}}
In this section, We work on minimal 3-folds of canonical index 1.
According to the Mori minimal model theory (\cite{KMM}, \cite{K-M}
and \cite{R1} etc.), one may take $X$ to be a minimal projective
threefold with invertible canonical sheaf $\omega_X={\mathcal
O}_X(K_X)$ and with locally factorial terminal singularities.
{}From the expression of the inequality in Theorem \ref{T1}. One
may assume $p_g(X)\ge 3$.

\begin{setup}\label{notations}{\bf Notations.} We study the canonical map $\varphi_{1}$
which is usually a rational map. Take the birational modification
$\pi: X'\rightarrow X$, according to Hironaka, such that

(i) $X'$ is smooth;

(ii) the movable part of $|K_{X'}|$ is basepoint free.

(iii) $\pi^*(K_X)$ is linearly equivalent to a divisor supported
by a divisor of normal crossings.

Denote by $g$ the composition $\varphi_{1}\circ\pi$. So $g:
X'\longrightarrow W'\subseteq{\Bbb P}^{p_g(X)-1}$ is a morphism.
Let $g: X'\overset{f}\longrightarrow B\overset{s}\longrightarrow
W'$ be the Stein factorization of $g$. We have the following
commutative diagram:

$$\begin{CD}
X' @>{f}>> B\\
@V{\pi}VV   @VV{s}V\\
X  @>>\varphi_{1}> W'
\end{CD}$$

We may write
$$K_{X'}=\pi^*(K_X)+\overline{E}=M+\overline{Z},$$
where $M$ is the movable part of $|K_{X'}|$, $\overline{Z}$ the
fixed part and $\overline{E}$ an effective divisor which is a sum
of distinct exceptional divisors. Throughout we always mean
$\pi^*(K_X)$ by $K_{X'}- \overline{E}$. So, whenever we take the
round up of $\alpha\pi^*(K_X)$, we always have
$\roundup{\alpha\pi^*(K_X)}\le \roundup{\alpha}K_{X'}$ for all
positive rational number $\alpha$. We may also write
$$\pi^*(K_X)=M+E_1',$$
where $E_1'=\overline{Z}-\overline{E}$ is actually an effective
divisor.

If $\dim\varphi_{1}(X)=2$, we see that a general fiber of $f$ is a
smooth projective curve of genus $g\ge 2$. We say that $X$ is {\it
canonically fibred by curves of genus $g$}.

If $\dim\varphi_{1}(X)=1$, we see that a general fiber $S$ of $f$
is a smooth projective surface of general type. We say that $X$ is
{\it canonically fibred by surfaces} with invariants $(c_1^2(S_0),
p_g(S)),$ where $S_0$ is the minimal model of $S$.
\end{setup}

\begin{setup}\label{B3}{\bf The case $\dim(B)=3.$} One has already the inequality
$$K_X^3\ge 2p_g(X)-6$$
according to Kobyashi (\cite{Kob}) even for a general minimal
3-fold of general type.
\end{setup}

For reader's convenience, we reformulate our known results in
\cite{CM} in the case $\dim (B)\le 2$.

\begin{thm}\label{known1} (Theorem 4.1 of \cite{CM})
Let $X$ be a minimal projective Gorenstein 3-fold of general type
with only locally factorial terminal singularities. Then we have

(i) If $\dim\varphi_{1}(X)=2$, {\it i.e.}, $X$ is canonically
fibred by curves of genus $g$,  then
$$K_X^3\ge \roundup{\frac{2}{3}(g-1)}(p_g(X)-2).$$

(ii) If $\dim\varphi_{1}(X)=1$, then either $K_X^3\ge 2p_g(X)-4$
or $(K_{S_0}^2, p_g(S))=(1, 2).$
\end{thm}

\begin{thm}\label{known2} (Theorem 4.3 of \cite{CM})
Let $X$ be a minimal projective smooth 3-fold of general type.
Suppose $\dim\varphi_{1}(X)=2$ and $X$ is canonically fibred by
curves of genus $2$. Then
$$K_X^3\ge \frac{2}{3}(2p_g(X)-5).$$
The inequality is sharp.
\end{thm}

We are left to study the only case when $\dim(B)=1$ and $X$ is
canonically fibred by surfaces with $(c_1^2, p_g)=(1,2)$. For this
purpose, we need a little bit of preparation.

\begin{setup}\label{main}{\bf Bounding $K_X^3$ from below.}
For the technical reason, we must assume that $K_X$ is ample from
now on. The inequality in Theorem \ref{T1} is trivial for small
value of $p_g(X)$. One may assume $p_g(X)\ge 3$. Furthermore, we
assume $\dim(B)=1$ and that a general fiber of the induced
fibration $f:X'\longrightarrow B$ is a surface with $(c_1^2,
p_g)=(1,2)$. Set $b:=g(B)$ the geometric genus of $B$.

By Lemma 4.5 of \cite{CM}, we have two cases exactly:
$$q(X)=b=1\ \ \text{and}\ \ h^2({\mathcal O}_X)=0,$$
$$q(X)=b=0\ \ \text{and}\ \ h^2({\mathcal O}_X)\le 1.$$
Write
$$|K_X|=|N|+Z$$
where $Z$ is the fixed part and $N$ the movable one. Then it is
obvious that
$$N=\pi_*(M)\ \ \text{and}\ \ Z=\pi_*(E_1').$$
Set $F:=\pi_*(S)$. One may write
$$M=\sum_{i=1}^aS_i$$
as a disjoint union of distinct fibers of $f$, where $a=p_g(X)-1$
whenever $b=0$, or $a=p_g(X)$ otherwise. Thus we have
$$N=\sum_{i=1}^aF_i$$
where $F_i=\pi_*(S_i).$

If $|N|$ has base points, then $F^2>0$ as a 1-cycle. Thus
$K_X\cdot F^2\ge 2$ because it is an even number. Then it is
obvious that
$$K_X^3\ge 2p_g(X)-2.$$

Otherwise, $|N|$ is base point free. In this case, $F$  is a
nonsingular projective surface with ample $K_F$. Also since
$Z|_F\sim K_F$ and $K_F^2=1$, we see that $Z|_F$ is an irreducible
curve on the surface $F$. Because $f$ obviously factors through
$X$, we denote by $f_1$ the induced fibration $X\longrightarrow
B$. Denote by $C$ the curve $Z\cap F$. Because $C\sim K_F$, $C$
must be a curve with arithmetical genus 2. Thus $C$ must be one of
the following types:

a) $C$ is smooth;

b) $C$ is an elliptic curve with exactly one node;

c) $C$ is an elliptic curve with exactly one cusp of type
$x^2=y^3$;

d) $C$ is a rational curve with exactly 2 nodes;

e) $C$ is a rational curve with one node and one cusp of type
$x^2=y^3$;

f) $C$ is a rational curve with exactly 2 cusps of type $x^2=y^3$;

g) $C$ is a rational curve with only one cusp of type $x^2=y^5$.

We will see later that the singularities on $C$ have  strong
connections with the value of $K_X^3$.

We consider the linear system $|K_{X'}+\pi^*(K_X)|$ on $X'$. For a
general fiber $S$, denote by $\sigma: S\longrightarrow F$ the
natural contraction.

Now we fix some notations. Assume $m>0$ is an integer and $0\le
n<m$. Under the premise of $a\ge m$, we may write $a:=a_1m+s$
where $a_1>0$ is an integer and $0\le s<m$. One may find distinct
smooth fibers $\{S_k\}$ such that
$$M\sim S_0+\sum_{j=1}^{a-a_1n-1}S_j+\sum_{i=a-a_1n}^{a}S_i.$$

Suppose the following condition (*) is satisfied:

(*) there is a number $r\ge 3$ such that, for all $i$ with
$a-a_1n\le i\le a$,
$$h^0(S_i,
K_{S_i}+\roundup{(\pi^*(K_X)-\frac{na_1}{ma_1+s}\pi^*(Z))|_{S_i})}\ge
r.$$

Because
$$\pi^*(K_X)-\sum_{i=a-a_1n}^{a}S_i-\frac{na_1}{ma_1+s}\pi^*(Z)\equiv
(1-\frac{na_1}{ma_1+s})\pi^*(K_X)$$ is nef and big, one has,
according to the Kawamata-Viehweg vanishing theorem, the exact
sequence:

\begin{align*}
0&\longrightarrow
H^0(X',K_{X'}+\roundup{\pi^*(K_X)-\sum_{i=a-a_1n}^aS_i-\frac{na_1}{ma_1+s}\pi^*(Z)})\\
 &\longrightarrow H^0(X',
K_{X'}+\roundup{\pi^*(K_X)-\frac{na_1}{ma_1+s}\pi^*(Z)})\\
&\longrightarrow\oplus_{i=a-a_1n}^{a}
H^0(S_i,K_{S_i}+\roundup{\pi^*(K_X)-\frac{na_1}{ma_1+s}\pi^*(Z)}|_{S_i}).
\end{align*}

Similarly because
$$\pi^*(K_X)-\sum_{j=1}^{a-a_1n-1}S_j-\sum_{i=a-a_1n}^{a}S_i-\frac{a-1}{a}\pi^*(Z)\equiv
\frac{1}{a}\pi^*(K_X)$$ is nef and big, one has again the exact
sequence:
\begin{align*}
0&\longrightarrow H^0(X',
K_{X'}+\roundup{\pi^*(K_X)-\sum_{j=1}^{a-a_1n-1}S_j-\sum_{i=a-a_1n}^{a}S_i-
\frac{a-1}{a}\pi^*(Z)})\\
&\longrightarrow H^0(X',K_{X'}+\roundup{\pi^*(K_X)-\sum_{i=a-a_1n}^aS_i-\frac{a-1}{a}\pi^*(Z)})\longrightarrow\\
&\oplus_{j=1}^{a-a_1n-1}H^0(S_j,
K_{S_j}+\roundup{\pi^*(K_X)-\sum_{i=a-a_1n}^aS_i-\frac{a-1}{a}\pi^*(Z)}|_{S_j}).
\end{align*}
Noting that
\begin{align*}
&K_{S_j}+\roundup{\pi^*(K_X)-\sum_{i=a-a_1n}^aS_i-\frac{a-1}{a}\pi^*(Z)}|_{S_j}\\
\ge &K_{S_j}+\roundup{(\pi^*(K_X)-\sum_{i=a-a_1n}^aS_i-\frac{a-1}{a}\pi^*(Z))|_{S_j}}\\
=&K_{S_j}+\roundup{\frac{1}{a}\pi^*(Z)|_{S_j}}\ge K_{S_j},\\
\end{align*}
we have
$$r_j:=h^0(S_j, K_{S_j}+\roundup{\frac{1}{a}\pi^*(Z)|_{S_j}})\ge
p_g(S_j)=2.$$

On the other hand, we have
\begin{align*}
h^0(X',
&K_{X'}+\roundup{\pi^*(K_X)-\sum_{j=1}^{a-a_1n-1}S_j-\sum_{i=a-a_1n}^{a}S_i-
\frac{a-1}{a}\pi^*(Z)})\\
\ge & h^0(X', K_{X'}+S_0).
\end{align*}

Whenever $h^2({\mathcal O}_X)=0$, the surjective map
$$H^0(X', K_{X'}+S_0)\longrightarrow H^0(S_0, K_{S_0})$$
gives
$$h^0(X', K_{X'}+S_0)=p_g(X)+2.$$
In this situation, we set $\delta:=2$.

Whenever $b=0$ and $h^2({\mathcal O}_X)=1$, we have
$$h^0(X', K_{X'}+S_0)\ge p_g(X)+1.$$
Whence we set $\delta:=1$.

The above two exact sequences give
\begin{equation*}
P_2(X)\ge p_g(X)+\delta+\sum_{j=1}^{a-a_1n-1}r_j+a_1nr.
\end{equation*}
Since
$$P_2(X)=\frac{1}{2}K_X^3-3(1-b+h^2({\mathcal O}_X)-p_g(X)),$$
one has
\begin{equation}\label{key}
K_X^3\ge
2(-2p_g(X)+\delta+\sum_{i=1}^{a-a_1n-1}r_j+a_1nr+3h^2({\mathcal
O}_X)-3b+3).
\end{equation}
\end{setup}

The above inequality is a key to better inequalities provided we
know all the numbers $r_j$, $r$ and $n$. We study it case by case
as follows. We first present the following

\begin{lem}\label{reduce}
On the general fiber $S$ of $f$, denote by
$D:=({\pi^*(Z)|_S})_{\text{red}}.$  Then $h^0(S, K_S+D)=2$ if and
only if $D$ is supported on a rational tree.
\end{lem}
\begin{proof}
Let $N$ be a very big natural number such that
$\roundup{\frac{1}{N}\pi^*(Z)|_S}=D$. Because
$$\pi^*(Z)|_S\sim \pi^*(K_X)|_S\cong \sigma^*(K_F)$$
is nef and big. Using the Riemann-Roch and the vanishing theorem,
we get
$$h^0(S, K_S+D)=\frac{1}{2}D\cdot (K_S+D)+\chi({\mathcal O}_S).$$
So $h^0(S, K_S+D)=2$ if and only if
$$(K_S+D)\cdot D=-2.$$
Because $D$ is 1-connected and reduced, it is the obvious fact
that $D$ supports on a rational tree. We are done.
\end{proof}

\begin{setup}\label{abc}{\bf The case \ref{main} a), b) and c).}
If there is a smooth fiber $F$ on $X$ such that $C=D\cap F$ is in
the case \ref{main} a), b) and c). Taking a smooth modification to
the morphism $f|_D:D\mapsto B$. One may easily see that $C$ is
always among these 3 types for a general fiber $F$. Thus we have
$r_j\ge 3$ for all $j$ by Lemma \ref{reduce}. We may take $n=0$.

Now if $b=1$, then the inequality (\ref{key}) gives
$$K_X^3\ge 2(p_g(X)-1).$$

If $b=0$ and $h^2({\mathcal O}_X)=1$, then (\ref{key}) gives
$$K_X^3\ge 2(p_g(X)+1).$$

If $b=0$ and $h^2({\mathcal O}_X)=0$, then (\ref{key}) gives
$$K_X^3\ge 2(p_g(X)-1).$$

We are left the situation that, for a general fiber $F$ of $f_1$,
$C$ falls into the cases \ref{main} d) through g). For these
cases, our argument depend on a special modification $\pi$.
\end{setup}

\begin{setup}{\bf The rest cases.}

{}From now on, we may suppose that $C$ is a singular rational
curve for a general fiber $F$ of $f_1$. We proceed our proof by
considering the singularities on the surface $Z$. First of all,
$Z$ must be singular along a curve. Otherwise, if $Z$ has isolated
singularities, $C$ would be a smooth curve of genus 2 which
contradicts to our assumption.

We hope to find a special embedded resolution of the pair $(X,Z)$
to prove Theorem \ref{T1}.

\begin{claim}\label{C1} $Z$ has at most 2 horizontal (with respect to $f_1$)
irreducible singular curves and the multiplicity of any such
singular curve on $Z$ is 2.
\end{claim}
\begin{proof}
In the process of finding the embedded resolution for $(X,Z)$, we
do not care those vertical modifications supported only on finite
number of fibers with regard to the fibration
$f_1:X\longrightarrow B$. This is because those vertical
modifications do not affect the behavior of $\pi^*(K_X)|_S$ for a
general fiber $S$ of $f$. By abuse of concepts, we call this kind
of vertical modifications to be {\it negligible}.

Pick up any irreducible singular curve $G$ of $Z$ such that
$f_1(G)=B$. Because $G$ has at most finite number of singular
points, we may take a negligible modification
$\pi_0:X_0\longrightarrow X$ such that $G$ is smooth upstairs.
Denote by $Z_0$ the strict transform of $Z$. We still denote by
$G$ the strict transform of $G$ upstairs. On $X_0$, $Z_0$ has a
singular curve along $G$ and $G$ is a smooth curve.

Let $\pi_1:X_1\longrightarrow X_0$ be a blow-up along the curve
$G$. Denote by $E_1$ the exceptional divisor on $X_1$. One may
write
$$\pi^*(Z_0)=Z_1+mE_1$$
where $Z_1$ is the strict transform of $Z_0$ and $m\ge 2$, because
$G$ belongs to singular locus of $Z_0$.

We consider the following commutative diagram:

$$\begin{CD}
S_{(1)} @>{\text{inclusion}}>> X_1\\
@V{\sigma_1}VV   @VV{\pi_0}V\\
S_{(0)}  @>>{\text{inclusion}}> X_0
\end{CD}$$
where $S_{(0)}=\pi_0^{-1}(F)$ and $S_{(1)}=\pi_1^{-1}(S_{(0)}).$
Denote by $Z_1\subset X_1$ the strict transform of $Z_0$. Then one
sees that $Z_1\cap S_{(1)}$ is irreducible. Also $\sigma_1$ is the
blow-up along the center $\{G\cap S_{(0)}\}$. Now we have
\begin{align*}
\sigma_{1}^*(C)&=\sigma_1^*(Z\cap S_{(0)})=\pi_0^*(Z)|_{S_{(1)}}\\
&=Z_1|_{S_{(1)}}+mE_1|_{S_{(1)}}.
\end{align*}
Because $m\ge 2$, we see that $G$ actually passes through a
singular point of $C$. Since $C$ has at most double points, $m\le
2$. Thus $m=2$ and $E_1|_{S_{(1)}}$ is either an irreducible
$(-1)$-curve or a sum of two distinct $(-1)$-curves. This also
means that, on $X_0$, one has $1\le G\cdot S_{(0)}\le 2$.

Now it is clear that $Z$ has at most two dictinct horizontal
singular curves like $G$ and the multiplicity of each is 2. The
lemma is proved.
\end{proof}

Based on the above argument, it is actually clear for us to
illustrate all possibilities.  Explicitly we have the following
possibilities:

I) if $Z_1$ is still singular along certain curve over $G$, then
$C_1$ is still singular and $C$ must have only one cusp (of type
$x^2=y^5$). In this case, $G$ is the only singular curve of $Z$
and $Z_1$ is singular along only one curve;

II) if $Z_1$ is smooth at generic points of $Z_1\cap E_1$ and the
natural map $\{E_1|_{Z_1}\}_{red}\mapsto G$ is not birational,
then $C$ has a node at each point of $\{G\cap F\}$.

III) if $Z_1$ is smooth at generic points of $Z_1\cap E_1$ and the
natural map $\{E_1|_{Z_1}\}_{red}\mapsto G$ is birational, then
$C$ has a cusp (of type $x^2=y^3$) at each point $\{G\cap F\}$.

\begin{lem}\label{node} For Case II), one has the same
inequalities as in \ref{abc}.
\end{lem}
\begin{proof}
If we are at Case II), then $\pi^*(K_X)|S=\pi^*(Z)|_S$ always
contains an elliptic cycle $C_0+C_1$ with $g(C_0)=g(C_1)=0$ and
$C_0\cdot C_1=2$ for a general fiber $S$. This means, for any $j$,
one has $r_j\ge 3$ by Lemma \ref{reduce}. We then take $n=0$ and
get the same inequalities as in \ref{abc}.
\end{proof}

Combining all arguments above, we are left the following 3
situations derived from possibilities I) and III):

A) $Z$ has only one horizontal singular curve $G$, $C$ has only
one cusp (of type $x^2=y^5$) for a general fiber $F$ and $G$ meets
the singular point of $C$;

B) $Z$ has only one horizontal singular curve $G$, $C$ has exactly
2 cusps (of type $x^2=y^3$) for a general fiber $F$, and $G$ meets
the 2 singular points of $C$;

B') $Z$ has two distinct horizontal singular curves $G$ and $H$,
$C$ has exactly 2 cusps (of type $x^2=y^3$) for a general fiber
$F$, and both $G$ and $H$ meet one singular point each of $C$.
\end{setup}

\begin{setup}\label{A}{\bf Embedded resolution of Type A).}
We construct an embedded resolution $\pi_{A}$ for the pair $(X,
Z)$ of Type A). Take $\pi_0$ and $\pi_1$ to be as in the proof of
Claim \ref{C1}. Because $Z_1$ still has a unique singular curve
which is over $G$, we denote such a curve by $G_1$. Modulo
negligible modifications, one may assume $G_1$ to be again
nonsingular.

Let $\pi_2: X_2\longrightarrow X_1 $ be the blow-up along $G_1$.
Denote by $Z_2$ the strict transform of $Z_1$ and by $E_2$ the
exceptional divisor. Set $S_{(2)}=\pi_2^{-1}(S_{(1)})$. Denote by
$\sigma_2:S_{(2)}\longrightarrow S_{(1)}$ the respective blow-up.
Because of the singularity type of $C$, one sees that the strict
transform $C_2$ of $C$ is already smooth. Since $E_2$ only touches
$C_2$ at one point, we denote by $G_2$ the reduced part of
$E_2|_{Z_2}$ which is of course irreducible. By considering the
multiplicity of $G_2$ in $Z_2$ (for instance taking blowing-ups
and then considering its impacts on $C_2$), one may see that $Z_2$
is smooth at generic points of $G_2$. Modulo negligible
modifications, $Z_2$ is already smooth. But the pull-back of $Z$
is in general not of normal crossing. We need more blow-ups.

Let $\pi_3:X_3\longrightarrow X_2$ be the blow-up along $G_2$
which, modulo negligible modifications, is a smooth curve. Let
$E_3$ be the exceptional divisor and $Z_3$ the strict transform of
$Z_2$. Denote by $S_{(3)}=\pi_3^{-1}(S_{(2)})$. Then one may see
that $E_2,$ $E_3$ and $Z_3$ meet at an irreducible curve $G_3$.

Finally blow-up $X_3$ along the curve $G_3$ (which could be smooth
modulo negligible modifications), we get $\pi_4:X_4\longrightarrow
X_3$. Denote by $E_4$ the exceptional divisor. Take more
negligible modifications, we get a resolution
$\pi_A:X'\longrightarrow X$ which is the composition of $\pi_i$
and those necessary negligible modifications. We replace our
original $\pi$ by $\pi_A$, keeping the same notations as above.
Pick up a general fiber $S$, then we may see that
$$\pi^*(K_X)|_S=\pi^*(Z)|_S=10L_4+5L_3+4L_2+2L_1+\tilde{C}$$
where $L_i$ are respective exceptional divisors induced from those
blow-ups of $\pi_i$ for $i=1,\ 2,\ 3,\ 4$ and $\tilde{C}$ is the
strict transform of $C$.

Denote by $\sigma:S\longrightarrow F$ the induced blow-up. Then
$\sigma^*(C)=\pi^*(Z)|_S$. {}From the whole process of blow-ups,
one sees that $\sigma^*(C)$ is a normal crossing divisor on $S$.
The intersection graph of $\sigma^*(C)$ is as follows:

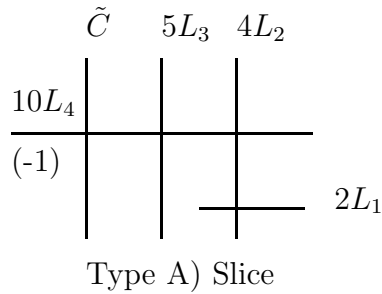
\setlength{\unitlength}{1mm}
\begin{picture}(70,40)
\thicklines \put(20,20){\line(1,0){40}}
\put(30,30){\line(0,-1){24}} \put(40,30){\line(0,-1){24}}
\put(50,30){\line(0,-1){24}} \put(45,10){\line(1,0){14}}
\put(20,23){$10L_4$} \put(30,33){$\tilde{C}$} \put(40,33){$5L_3$}
\put(50,33){$4L_2$} \put(63,10){$2L_1$} \put(20,15){(-1)}
\put(30,0){Type A) Slice}
\end{picture}
\medskip

Take two joint objects from $\{\tilde{Z},\ E_1,\ E_2,\ E_3,\
E_4\}$, they meet a general fiber $S$ at exactly one point.
 According to the next lemma, the divisor
 $$\tilde{Z}+E_1+E_2+E_3+E_4$$
 is normal crossing over a general point of $B$. Taking necessary  negligible
 modifications, $\pi_A$ is finally an embedded resolution of $(X,Z)$.
\end{setup}

\begin{lem}\label{nc} Let $P$ be a point of a smooth variety $V$.
Suppose I have 3 irreducible smooth divisor $H_i\subset V$ for
$i=1,\ 2,\ 3$ such that $P\in H_i$ for all $i$. Assume
$$(H_1\cdot H_2\cdot H_3)_P=1.$$
Then $H_1+H_2+H_3$ is a normal crossing divisor at $P$.
\end{lem}
\begin{proof}
This is a trivial statement. Denote by $f_i$ the local equations
of $H_i$ for all $i$. Then, by definition,
$$\dim_{k}({\mathcal O}_{V,P}/(f_1, f_2, f_3))=(H_1\cdot H_2\cdot
H_3)_P=1.$$ This means that $f_1$, $f_2$ and $f_3$ actually form a
local parameters of the point $P$. We are done.
\end{proof}

\begin{setup}\label{B}{\bf Embedded resolution of Type $\text{B}^*$).}
The construction is somehow similar to \ref{A}. Both Type B) and
Type B') are essentially the same case. We omit those minor
differences of the details for Type B) which is simply a copy of
the one below. The most important point is that, for Type B) and
Type B'), we finally have the same $\pi^*(Z)|_S$ for a general
fiber $S$.

Take $\pi_0$ to be as in the proof of Claim \ref{C1}. Let $\pi_1:
X_1\longrightarrow X_0$ be the blow-up along two smooth curve $G$
and $H$. Denote by $Z_1$ the strict transform of $Z$, and by
$E_1$, $E_1'$ the exceptional divisors. Similarly, one may see
that $Z_1$ is already smooth simply because of the singularity
type of $C$. We also see that $E_1$ (or $E_1'$) and $Z_1$ meet at
an irreducible curve $G_1$ (or $H_1$). Modulo negligible
modifications, one may assume both $G_1$ and $H_1$ are smooth
curves. We keep parallel notations as in \ref{A}.

Going on blow-ups along $G_1$ and $H_1$, one gets
$\pi_2:X_2\longrightarrow X_1$. Denote by $E_2$, $E_2'$ the
exceptional divisors. One sees that $E_2$ (or $E_2'$), $E_1$ (or
$E_1'$) and $Z_2$ still meet along an irreducible curve $G_2$ (or
$H_2$). One may take negligible modifications such that $G_2$ and
$H_2$ are smooth.

We finally blow-up $X_2$ along $G_2$ and $H_2$ to obtain
$\pi_3:X_3\longrightarrow X_2$. Taking further negligible
modification, we get a resolution $\pi_B: X'\longrightarrow X$. We
have
$$\pi_B^*(Z)=\tilde{Z}+6E_3+6E_3'+3E_2+3E_2'+2E_1+2E_1'.$$
The slice on a general fiber $S$ is
$$\pi_B^*(Z)|_S=\tilde{C}+6L_3+6L_3'+3L_2+3L_2'+2L_1+2L_1'$$
where $L_i=E_i|_S$ and $L_i'=E_i'|_S$ for all $i$. The
intersection graph is as follows.

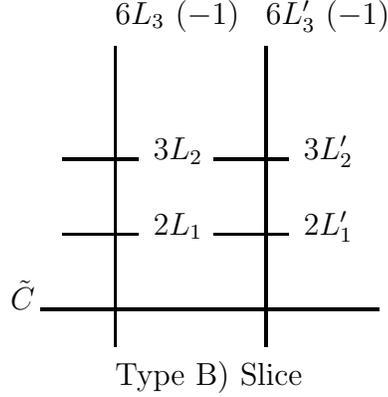
\setlength{\unitlength}{1mm}
\begin{picture}(80,50) \thicklines
\put(30,10){\line(1,0){45}} \put(26,10){$\tilde{C}$}
\put(40,5){\line(0,1){40}} \put(40,48){$6L_3\ (-1)$}
\put(60,5){\line(0,1){40}} \put(60,48){$6L_3'\ (-1)$}
\put(33,20){\line(1,0){10}} \put(45,20){$2L_1$}
\put(53,20){\line(1,0){10}} \put(65,20){$2L_1'$}
\put(33,30){\line(1,0){10}} \put(45,30){$3L_2$}
\put(53,30){\line(1,0){10}} \put(65,30){$3L_2'$} \put(40,0){Type
B) Slice}
\end{picture}

Applying Lemma \ref{nc}, one may see that $\pi_B$ is an embedded
resolution of $(X,Z)$.
\end{setup}

\begin{setup}{\bf The inequalities for Type A) case.}
We apply the argument in \ref{main}. We take two integers $m$ and
$n$ such that $1-\frac{na_1}{ma_1+s}>\frac{1}{10}.$ Write
$a:=a_1m+s$ as in \ref{main}. In order to get an effective
inequality, we only need to verify the condition (*). For
simplicity, we still denote $\pi_A$ by $\pi$. Recall that we have
$$\pi^*(K_X)|_S=\pi^*(Z)|_S=10L_4+5L_3+4L_2+2L_1+\tilde{C}$$
where $L_4^2=-1$ and all these curves are smooth rational curves.
Set $D_0:=(\pi^*(Z)|_S)_{red}$. Then $D_0$ is of course a rational
tree. {}From the intersection form of $\pi^*(Z)|_S$, we have
$$D_0\cdot L_4=2.$$
Now we verify the condition (*). For a general fiber $S$, we have
\begin{align*}
&h^0(S, K_S+\roundup{(\pi^*(K_X)-\frac{na_1}{ma_1+s}\pi^*(Z))|_S})\\
\ge&h^0(S, K_S+D_0+L_4)=h^0(S, K_S+\roundup{\frac{1}{5}\pi^*(Z)|_S})\\
=&\frac{1}{2}(K_S+D_0+L_4)(D_0+L_4)+\chi({\mathcal O}_S)=3.\\
\end{align*}

Thus the inequality (\ref{key}) gives
\begin{equation}\label{TA}
K_X^3\ge 2(-2p_g(X)+\delta+(2+\frac{n}{m})a+3h^2({\mathcal
O}_X)-3b-\frac{n(m-1)}{m}+1).
\end{equation}

Now we see that, if $p_g(X)$ is very big (thus $m$ can be big),
then the ratio $K_X^3/p_g(X)$ is close to $\frac{9}{5}$, a very
good inequality. The inequality in Theorem \ref{T1} allows us to
assume $p_g(X)\ge 5$ and so $a\ge 4$. We may take $m=4$ and $n=3$.
Apparently, $1-\frac{3a_1}{4a_1+s}\ge \frac{1}{7}>\frac{1}{10}$
because $s<4$ by definition.

Explicitly, if $b=1$, then we have
$$K_X^3\ge \frac{3}{2}(p_g(X)-3).$$
This is better than the inequality in Theorem \ref{T1} only
whenever $p_g(X)\ge 7$. But the trivial inequality $K_X^3\ge
p_g(X)$ amends whenever $p_g(X)\le 6$.

If $b=0$ and $h^2({\mathcal O}_X)=0$, then we have
$$K_X^3\ge 2(\frac{3}{4}p_g(X)-2).$$
This is better than the inequality in Theorem \ref{T1}.

If $b=0$ and $h^2({\mathcal O}_X)=1$, then we have
$$K_X^3\ge \frac{3}{2}p_g(X).$$
This is much better than what we want in Theorem \ref{T1}.
\end{setup}

\begin{setup}{\bf The inequalities for Type B), Type B') case.}
We take positive integers $m$ and $n$ such that
$1-\frac{na_1}{ma_1+s}>\frac{1}{6}$. Write $a:=a_1m+s$ as in
\ref{main}. In order to get an effective inequality, we only need
to verify the condition (*). For simplicity, we still denote
$\pi_B$ by $\pi$. Recall that we have
$$\pi^*(K_X)|_S=\pi^*(Z)|_S=6(L_3+L_3')+3(L_2+L_2')+2(L_1+L_1')+\tilde{C}$$
where $L_3^2={L_3'}^2=-1$ and all these curves are smooth rational
curves. Set $D_0:=(\pi^*(Z)|_S)_{red}$. Then $D_0$ is of course a
rational tree. {}From the intersection form of $\pi^*(Z)|_S$, we
have
$$D_0\cdot L_3=D_0\cdot L_3'=2.$$
Now we verify the condition (*). For a general fiber $S$, we have
\begin{align*}
&h^0(S, K_S+\roundup{(\pi^*(K_X)-\frac{na_1}{ma_1+s}\pi^*(Z))|_S})\\
\ge &h^0(S, K_S+D_0+L_3+L_3')=h^0(S, K_S+\roundup{\frac{1}{3}\pi^*(Z)|_S})\\
=&\frac{1}{2}(K_S+D_0+L_3+L_3')(D_0+L_3+L_3')+\chi({\mathcal O}_S)=4.\\
\end{align*}

Thus the inequality (\ref{key}) gives
\begin{equation}\label{TB}
K_X^3\ge 2(-2p_g(X)+\delta+(2+\frac{2n}{m})a+3h^2({\mathcal
O}_X)-3b-\frac{2n(m-1)}{m}+1).
\end{equation}

Still one may see that, if $p_g(X)$ is bigger, the ratio
$K_X^3/p_g(X)$ is close to $\frac{10}{3}$. Under the assumption of
$p_g(X)\ge 4$, we may take $m=3$ and $n=2$. So
$$1-\frac{2a_1}{3a_1+s}\ge \frac{1}{5}>\frac{1}{6}.$$

Explicitly, if $b=1$, then we have
$$K_X^3\ge \frac{8}{3}(p_g(X)-2).$$
This is better than what we want in Theorem \ref{T1}.

If $b=0$ and $h^2({\mathcal O}_X)=0$, then we have
$$K_X^3\ge 2(\frac{4}{3}p_g(X)-3).$$
This is also better than what we want in Theorem \ref{T1}.

If $b=0$ and $h^2({\mathcal O}_X)=1$, then we have
$$K_X^3\ge \frac{8}{3}p_g(X)-2.$$

Already these inequalities are better than the one in Theorem
\ref{T1}.
\end{setup}

\begin{setup}{\bf Summary.}
Comparing what we have got, we may conclude Theorem \ref{T1}.
\end{setup}

\begin{setup}{\bf Proof of Corollary \ref{cor}.}
\begin{proof}
If $K_X^3<\frac{3}{2}p_g(X)-\frac{9}{2}$, then we have $p_g(X)\ge
5$. This means that we have the canonical map. {}From the argument
above, the only possibility is that $X$ is canonically fibred by
curves of genus $2$.  This inequality is not empty according to
Kobayashi's example.
\end{proof}
\end{setup}

\begin{setup}{\bf Proof of Theorem \ref{T1.5}}
\begin{proof}
If $K_X^3<2p_g(X)-6$, then one has $p_g(X)\ge 5$. We may study the
canonical map. Both \ref{B3} and Theorem \ref{known1} tell that
either $X$ is canonically fibred by curves of genus $2$ or $X$ is
canonically fibred by surfaces of general type with $(c_1^2,
p_g)=(1,2)$. We study the later case.

We take the induced fibration $f: X'\longrightarrow  B$ where a
general fiber $S$ is a smooth projective surface with $(c_1^2,
p_g)=(1,2)$. Noting that $f_*\omega_{X'/B}$ is a vector bundle of
rank $2$ because $p_g(S)=2$, we considering the natural projection
$$p: \bP(f_*\omega_{X'/B})\longrightarrow B.$$
Because $K_{X'}+S_1+S_2\ge K_{X'}$, we see that the fibration $f$
rationally factors through $p$.  Taking birational modifications,
we may have a morphism from $f$ to $p$. Thus we have the following
commutative diagram:

$$\begin{CD}
X' @>{\psi}>> \bP(f_*\omega_{X'/B})\\
@V{f}VV   @VV{p}V\\
B  @>>{\text{identity}}> B
\end{CD}$$

For any fiber $S$ of $f$, we see that $\psi|_S=\phi_{K_S}$ because
we have
$$|K_{X'}+S_1+S_2||_S=|K_S|$$
by Lemma 4.6 of \cite{CM}, where the $S_i$ are general smooth
fibers of $f$. Therefore $\psi$ is actually a fibration over a
ruled surface with a general fiber a smooth curve of genus $2$. We
are done.
\end{proof}
\end{setup}

\section{\bf Proof of Theorem \ref{T2}}

Though partial effective Noether type of inequalities for a
general minimal 3-fold are given in \cite{CM}, there remain
several hard cases to study. In this section, we are able to
develop the technique in \cite{IJM} to present integral and more
precise results which make it possible for us to describe those
3-folds with small slope $\frac{K_X^3}{p_g(X)}$.

\begin{setup}\label{symbols}{\bf Notations.} We are treating a
general object so that most of the divisors we come across are
rational divisors. In order to prove Theorem \ref{T2}, we may
assume that $X$ is a normal projective minimal 3-fold with only
$\bQ$-factorial terminal singularities. We suppose $p_g(X)\ge 2$.

We study the canonical map $\varphi_{1}$ which is usually a
rational map. Take the birational modification $\pi: X'\rightarrow
X$, according to Hironaka, such that

(i) $X'$ is smooth;

(ii) the movable part of $|K_{X'}|$ is basepoint free. (Sometimes
we even call for such a modification that those movable parts of a
finite number of linear systems are all basepoint free.)

(iii) $\pi^*(K_X)$ is linearly equivalent to a divisor supported
by a divisor of normal crossings.

Denote by $g$ the composition $\varphi_{1}\circ\pi$. So $g:
X'\longrightarrow W'\subseteq{\Bbb P}^{p_g(X)-1}$ is a morphism.
Let $g: X'\overset{f}\longrightarrow B\overset{s}\longrightarrow
W'$ be the Stein factorization of $g$. So we have the same
commutative diagram as in \ref{notations}. Write
$$K_{X'}=_{\Bbb Q}\pi^*(K_X)+E_1=_{\Bbb Q}M_1+Z_1,$$
where $M_1$ is the movable part of $|K_{X'}|$, $Z_1$ the fixed
part and $E_1$ an effective ${\Bbb Q}$-divisor which is a ${\Bbb
Q}$-sum of distinct exceptional divisors. Throughout we always
mean $\pi^*(K_X)$ by $K_{X'}- E_1$. So, whenever we take the round
up of $m\pi^*(K_X)$, we always have $\roundup{m\pi^*(K_X)}\le
mK_{X'}$ for all positive number $m$. We may also write
$$\pi^*(K_X)=_{\Bbb Q} M_1+E_1',$$
where $E_1'=Z_1-E_1$ is actually an effective ${\Bbb Q}$-divisor.

If $\dim\varphi_{1}(X)=2$, we see that a general fiber of $f$ is a
smooth projective curve $C$ of genus $g\ge 2$. If
$\dim\varphi_{1}(X)=1$, we see that a general fiber $S$ of $f$ is
a smooth projective surface  $S$ of general type. The invariants
of $S$ are  $(c_1^2(S_0), p_g(S))$
 where $S_0$ is the minimal model of $S$.

{\it A generic irreducible element $S$ of} $|M_1|$ means either a
general member of $|M_1|$ whenever $\dim\varphi_{1}(X)\ge 2$ or,
otherwise, a general fiber of $f$.
\end{setup}

For reader's convenience, we recall known results from \cite{CM}.
\begin{thm}\label{known3} (Theorem 3 and Proposition 3.4 of
\cite{CM}) Under the above assumptions, one has

1) $K_X^3\ge 2p_g(X)-4$ whenever $p_g(X)\ge 6$, $\dim(B)=2$ and
$g(C)\ge 3$;

2) $K_X^3\ge \frac{3}{2}p_g(X)-\frac{5}{2}$ whenever $p_g(X)\ge
10$, $(c_1^2(S_0), p_g(S))=(1,1)$, $\dim(B)=1$ and
$\dim\varphi_{2K_X}(X)\ge 2$.
\end{thm}

Standard surface theory tells us that a surface $S$ of general
type with $K_{S_0}^2=1$ has only 2 possibilities: either
$p_g(S)=1$ or $p_g(S)=2$.

Before proving Theorem \ref{T2}, we must study the other cases.
The following proposition presents a general method to estimate
certain intersection numbers on $X$.

\begin{prop}\label{estimate}
Let $X$ be a minimal projective 3-fold of general type with only
${\Bbb Q}$-factorial terminal singularities and assume $p_g(X)\ge
2$. Keep the same notations as in \ref{symbols} . Pick up a
generic irreducible element $S$ of $|M_1|$. Suppose, on the smooth
surface $S$, there is a movable linear system $|G|$ and denote by
$C$ a generic irreducible element of $|G|$. Set
$\xi:=(\pi^*(K_X)\cdot C)_{X'}$ and
$$p:=\begin{cases} 1  &\text{if}\ \dim\varphi_1(X)\ge 2\\
a  &\text{if}\ \ \pi^*(K_X)\equiv_{\bQ}aS+\text{effective}\ \
\bQ\text{-divisors}
\end{cases}$$

Assume

(i) there is a rational number $\beta>0$ such that
$\pi^*(K_X)|_S-\beta C$ is numerically equivalent to an effective
${\Bbb Q}$-divisor;

(ii) the inequality
$\alpha:=(m-1-\frac{1}{p}-\frac{1}{\beta})\xi>1$ holds. Set
$\alpha_0:=\roundup{\alpha}$. Then we have the inequality
$$m\xi\ge 2g(C)-2+\alpha_0$$.
\end{prop}

\begin{proof} This is a weak version of Theorem 2.2 in \cite{IJM}.
We do not need the birationality of $\varphi_m$. So one may drop
additional assumptions there.
\end{proof}

\begin{setup}\label{K>2}{\bf The case $\dim(B)=1$ and
$c_1^2(S_0)\ge 2$.} We have
$$
\pi^*(K_X)=_{\bQ} M_1+E_1' \equiv_{\bQ}aS+E_1'$$ where $a\ge
p_g(X)-1$. So one has
$$K_X^3=\pi^*(K_X)^3\ge(\pi^*(K_X)^2\cdot S)(p_g(X)-1).$$

If $b=g(B)>0$, then the movable part of $|K_X|$ is already base
point free. Thus one has
$$\pi^*(K_X)|_S=\sigma^*(K_{S_0}).$$
Thus $\pi^*(K_X)^2\cdot S=(\sigma^*(K_{S_0}))^2\ge 2$. So
$$K_X^3\ge 2(p_g(X)-1).$$

{}From now on, we assume $b=0$ and $p_g(X)\ge 12$. In order to
apply Proposition \ref{estimate}, we must find the number $\beta$
and the curve $C$.

Note that $p_g(X)>0$ implies $p_g(S)>0$. According to \cite{Ci},
we know that $|2K_{S_0}|$ is base point free. So is
$|2\sigma^*(K_{S_0})|$. We set $C$ be a general member of
$|2\sigma^*(K_{S_0})|$. So $C$ is a smooth curve with
$\deg(K_C)\ge 12$. According to Step 2 of Proposition 3.3 in
\cite{CM}, we have
$$\pi^*(K_X)|_S\ge_{\bQ}\frac{5}{6}\sigma^*(K_{S_0}).$$
Thus we may set $\beta=\frac{5}{12}$. Also one may set $p=11$. An
initial lower bound for $\xi$ is
$$\xi\ge \frac{5}{12}C^2\ge \frac{10}{3}.$$
Now we may choose $m$ and run Proposition \ref{estimate}.

Take $m_1=4$. Then $\alpha_1\ge
(3-\frac{1}{11}-\frac{12}{5})\xi\ge \frac{64}{33}$. So
$\alpha_0\ge 2$. Applying Proposition \ref{estimate}, one gets
$\xi\ge \frac{7}{2}$.

Take $m_2=5$. Then $\alpha_2=(4-\frac{1}{11}-\frac{12}{5})\xi\ge
\frac{581}{110}$. So $\alpha_0\ge 6$. Applying Proposition
\ref{estimate}, one gets $\xi\ge \frac{18}{5}$.

Take $m_3=6$. Then $\alpha_3=(5-\frac{1}{11}-\frac{12}{5})\xi\ge
\frac{2484}{275}>9$. So $\alpha_0\ge 10$. Applying Proposition
\ref{estimate}, one gets $\xi\ge \frac{11}{3}$ which might be the
best bound through our method.

So we have $(\pi^*(K_X)|_S)^2\ge \frac{5}{12}\xi\ge
\frac{55}{36}>\frac{3}{2}.$ Thus we have the inequality
\begin{equation}\label{K2}
K_X^3\ge \frac{55}{36}(p_g(X)-1).\end{equation}
\end{setup}

\begin{setup}\label{K1}{\bf The case $\dim(B)=1$,
$K_{S_0}^2=p_g(S)=1$ and $\dim\varphi_{2k_X}(X)=1$.} We assume
$p_g(X)\ge 10$.  Considering the induced fibration
$f:X'\longrightarrow B$, we have $q(X)\le 1$ and
$q(X)-h^2({\mathcal O}_X)\ge 0$ according to \cite{JMSJ2}. In
fact, this case is very simple since $f_*\omega_{X'}$ is an
invertible sheaf while $R^1f_*\omega_{X'}=0$. So we have
$$\chi({\mathcal O}_X)=1-q(X)+h^2({\mathcal O}_X)-p_g(X)\le
1-p_g(X).$$ Applying Reid's plurigenus formula (\cite{YPG}), one
has
\begin{align*}P_2(X)&\ge \frac{1}{2}K_X^3-3\chi({\mathcal O}_X)\\
&\ge \frac{1}{2}K_X^3+3p_g(X)-3.
\end{align*}
We may remodify our original $\pi$ such that the movable part of
$|2K_{X'}|$ is also base point free. Write
$$|2K_{X'}|=|M_2|+Z_2$$
where $M_2$ is the movable part. One has
$$M_2\equiv a_2S$$
where $a_2\ge P_2(X)-1$.  Because
$$2\pi^*(K_X)\ge_{\bQ} M_2,$$
we have
\begin{equation}\label{bu} 2K_X^3\ge
a_2(\pi^*(K_X)|_S)^2.
\end{equation}
 So the key point might be to estimate the
number $(\pi^*(K_X)|_S)^2$ which is a rational number.

The base point freeness of $|2\sigma^*(K_{S_0})|$ allows us to
take $C$ to be a general member of this system. Then $C$ is a
smooth curve with $\deg(K_C)=6$. Because
$$\pi^*(K_X)\equiv_{\bQ}\frac{a_2}{2}S+*,$$
we may take $p=\rounddown{\frac{a_2}{2}}\ge 12.$ Similarly
$5\pi^*(K_X)|_S\ge_{\bQ}4\sigma^*(K_{S_0})$ by Step 2 of
Proposition 3.3 in \cite{CM}, we may take $\beta=\frac{2}{5}$. We
have $\xi\ge \frac{2}{5}C^2\ge\frac{8}{5}.$

Now take $m_1=5$. Then
$$\alpha_1=(4-\frac{1}{p}-\frac{1}{\beta})\xi>2.$$
Proposition \ref{estimate} gives $\xi\ge \frac{9}{5}$.

Take $m_2=6$. Then
$$\alpha_2=(5-\frac{1}{p}-\frac{1}{\beta})\xi>4.$$
Proposition \ref{estimate} gives $\xi\ge \frac{11}{6}$.

Take $m_3=7$. Then
$$\alpha_3=(6-\frac{1}{p}-\frac{1}{\beta})\xi>6.$$
Proposition \ref{estimate} gives $\xi\ge \frac{13}{7}$.

In general, we may get
$$\xi\ge \frac{2m_k-1}{m_k}$$
for all $m_k\ge 8$ by induction. Thus $\xi\ge 2$. This means
$(\pi^*(K_X)|_S)^2\ge \frac{2}{5}\xi\ge \frac{4}{5}.$ So the
inequality (\ref{bu}) becomes
\begin{equation}\label{K=1}
K_X^3\ge \frac{3}{2}p_g(X)-2.
\end{equation}
\end{setup}

\begin{setup}{\bf Proof of Theorem \ref{T2}.}
\begin{proof}
Assume $K_X^3<\frac{3}{2}p_g(X)-\frac{5}{2}.$ Because $K_X^3>0$,
one sees that $p_g(X)\ge 2$. So one may always consider the
canonical map $\varphi_1$. Suppose $p_g(X)\ge 12$. According to
\ref{B3}, Theorem \ref{known3} and the inequalities (\ref{K2}) and
(\ref{K=1}), $X$ must be either canonically fibred by curves of
genus $2$ or canonically fibred by surfaces of general type with
$(c_1^2, p_g)=(1,2)$. Then a parallel argument to that in the
proof of Theorem \ref{T1.5} also works. We are done.
\end{proof}
\end{setup}

In fact, a combination of \cite{CM} and this section may present
the following more general result for which we omit the details.
\begin{thm}\label{T2'} There are two sequences (computable) of positive
rational numbers $\{a_k\}$ and $\{b_k\}$ with $k\ge 2$ such that

1) $\frac{4}{3}<a_k\le 2$ for all $k\ge 2$ and $a_{k_1}\le
a_{k_2}$ whenever $k_1<k_2$;

2) $\lim_{k\rightarrow +\infty}a_k=2$ and $\{b_k\}$ is bounded;

3) for any minimal projective 3-fold $X$ of general type with
canonical singularities, set $k:=\rounddown{\frac{p_g(X)-2}{2}}.$
If $K_X^3<a_kp_g(X)-b_k$ and $p_g(X)\not\in [2,5]$, then $X$ is
fibred by curves of genus 2.
\end{thm}

\begin{setup}{\bf Examples.} The only known examples satisfying the equality in
Theorem \ref{T1} or the assumption of both Theorem \ref{T1.5} and
Theorem \ref{T2} were found by M. Kobayashi (\cite{Kob}).
\end{setup}

\begin{setup}{\bf An open problem.}
If $X$ is a Gorenstein minimal projective threefold of general
type with only canonical singularities, then it is well-known that
$\chi({\mathcal O}_X)<0$ according to Miyaoka (\cite{M}). There
should be an analogue of the Noether inequality as in Theorem
\ref{T1} in the form:
$$K_X^3\ge -a\chi({\mathcal O}_X)-b$$
where $a$ and $b$ are positive rational numbers. One may try to
study the bicanonical map of $X$. We have an effective lower bound
for $a$. Any bound $a>1$ is nontrivial and interesting. The
author's opinion is that to find a Noether inequality in this
direction is more difficult simply because the inter relations
among $p_g$, $q$ and $h^2({\mathcal O}_X)$ are far from being
clear to us, unlike in surface case.
\end{setup}

%%%%%%%%%%%%%%%%%%%%%%%%%%%%%%%%%%%%%%%%%%%%%%%%%%%%%%%

\end{document}